\documentclass[A4paper,11pt]{article}

\pdfoutput=1

\usepackage{latexsym}
\usepackage{epsfig}
\usepackage{verbatim}
\usepackage{multirow}
\usepackage{fancybox}
\usepackage{shadow}
\usepackage{wrapfig, lipsum}
\usepackage{cite}
\usepackage{amssymb}
\usepackage{amsmath}
\usepackage{bm}
\usepackage{amscd}
\usepackage{graphicx}
\usepackage{xcolor}

\newtheorem{theorem}{\sc Theorem}[section]
\newtheorem{lemma}[theorem]{\sc Lemma}

\newcommand{\eps}{\varepsilon}

\newcommand{\Varr }{\mathrm{var}\,}
\newcommand{\Cov}{\mathrm{cov}\,}

\newcommand{\proofend}{{\medskip\medskip}}
\newcommand{\proof}{{\noindent\em Proof. }}

\author{
  {\sc Bernard Chazelle}
\thanks{Department of Computer Science,
       Princeton University, 
{\tt chazelle}@{\tt cs.princeton.edu }}
\and
{\sc Chu Wang}
\thanks{ Nokia Bell Labs,
{\tt   chu.wang}@{\tt nokia.com }}
}

\title{Self-Sustaining Iterated Learning
\thanks{This work was supported in part by NSF grant CCF-1420112.
}}

\vspace{.5cm}

\date{}

\begin{document} \maketitle
%\icmlkeywords{Iterated learning, language evolution, opinion dynamics, convergence, sustainable learning}

\vskip 0.3in

\begin{abstract}
An important result from psycholinguistics (Griffiths \& Kalish, 2005) states that no language can be learned
iteratively by rational agents in a self-sustaining manner.   
We show how to modify the learning process slightly 
in order to achieve self-sustainability.
Our work is in two parts. First, we characterize iterated learnability in geometric terms
and show how a slight, steady increase in the lengths of the training sessions
ensures self-sustainability for any discrete language class.
In the second part, we tackle the nondiscrete case and investigate
self-sustainability for iterated linear regression.
We discuss the implications of our findings to issues
of non-equilibrium dynamics in natural algorithms. 
\end{abstract}

\section{Introduction}

Consider this hypothetical scenario: A native
speaker of Quenya\footnote{Quenya
is one of J.R.R. Tolkien's fictional languages.}
sets out to teach the language to an English speaker; after 
a year of teaching, the learner considers herself fluent enough to teach Quenya
to some other English speaker, who a year later does the same.
In this form of {\em iterated learning}, agents teach each other in sequence:
X teaches Y, who then teaches Z, who then 
teaches...\cite{beppu2009iterated,griffiths2007language,griffiths2005bayesian,rafferty2009convergence,kirby2014iterated,griffiths2008theoretical,perfors2011language,rafferty2014analyzing,smith2009iterated,kalish2007iterated}.
By a classic result of Griffiths and Kalish~\cite{griffiths2005bayesian},
Quenya will vanish after a finite number of iterations,
at which point the agents, assumed to be rational, will be ``teaching"
each other plain English. In other words, after a while,
learners will be taught nothing they don't already know:
iterated learning is not self-sustaining.

Such findings are hard to validate empirically but variants of it are within the reach
of experimental psychology.
As early as 1932, in fact, the English psychologist Frederic Bartlett used iterated
learning to expose hidden biases among humans.  He presented a picture of an owl to
a person for given period of time and then asked her to draw it from memory.
Her picture was then shown to the next learner for the same amount of time,
who then proceeded to draw it back from memory. After 20 iterations of this process,
to Bartlett's surprise, 
what was being drawn was no longer an owl but, quite clearly, a cat! 
The challenge was to explain why humans would exhibit a pro-feline bias
without falling into the trap of just-so stories. 

Griffiths et al.~\cite{kalish2007iterated} repeated the same experiment ten years ago,
this time trading owls for lines. The goal was to see
if linear regression could be iterated: the answer was a resounding No. 
Skipping over logistical details, the experiment presents the first learner with
a cloud of 20 points drawn randomly, with noise, from the line $Y= 1-X$.
The cloud vanishes and the learner is then asked to reconstruct it from memory.
She then becomes the teacher by 
passing on her own cloud to the next learner, who likewise, looks at it for a while,
and then tries to reconstruct it from memory, etc. 
Surprisingly, iterating this process a mere nine times leads the last learner 
in the sequence to draw a cloud that regresses to the line $Y=X$; in other words, teaching
about descending lines iteratively has precisely the opposite effect!
In fact the initial picture is essentially is irrelevant. A random cloud of points will also
lead to $Y=X$. 

Unlike the Quenya scenario, where the bias toward English
is not unexpected, the cat and line experiments both 
reveal a hidden prior among the participants. 
Humans seem to love cats and possess a strong positive correlation bias;
it is easy to speculate why.\footnote{Our favorite piece of
anecdotal evidence in support of the positive slope bias is that  
no road sign in the US features an aircraft on a descending path.}
It is noteworthy 
that the prior should prevail even in the absence of any sort of priming.
Indeed, this experiment fails miserably if you try it yourself
by playing the role of all the agents in sequence. The use of 
different learners ensures that the training does not acquire long-term memory.
Similar laboratory experiments with human subjects (well, undergraduates) have confirmed the 
unstainability of iterated learning~\cite{kalish2007iterated,beppu2009iterated,tamariz2015culture,bartlettremembering,griffiths2008theoretical}.

In our first example, 
Quenya gets ``washed out" by English in a way reminiscent of the fixation of an allele
through genetic drift. Indeed, the original impetus 
for studying iterated learning in psycholinguistics was to look for 
a parallel to Kimura's neutral theory of molecular evolution
in the area of cultural transmission. 
People learn their native tongue from speakers who themselves learned it from others.
This process introduces variation along the way, some of which is retained
durably. The selectionist view seeks to explain this 
process by fitness considerations at the population level.
Iterated learning suggests a different explanation. Language acquisition
suffers from a well-documented information bottleneck (the notorious ``poverty of stimulus"),
so one might expect languages to evolve so as to be easy to learn:  could 
complexity theory be the key?  This push for simplicit would then 
trigger the emergence of linguistic universals (eg, compositionality) that one finds 
present in all languages~\cite{griffiths2007language}.
This view complements---some will argue, contradict---Chomsky's interpretation of universals as the product of
constraints imposed by an innate genetic endowment.

Following Chomsky and Lasnik's theory of ``Principles and Parameters,"
Rafferty et al.~\cite{rafferty2009convergence} model languages
by means of a handful of parameters: think of a few knobs whose settings specify 
any given language.  Language evolution thus entails the trans-generational update
of a probability distribution over that parameter space. Assuming that the learners are rational
Bayesian agents, iterated learning acts as a Gibbs sampler 
for a joint probability distribution over languages and their sentences. By converging to
a stationary distribution, iterated learning proves 
incapable of sustaining itself past the mixing time.   In that model, languages
evolve to reflect the priors of the learners while losing all trace of the ancestor language.
While this phenomenon is of central relevance in the study of universal grammars,
it leaves open the possibility that changes in the sampling algorithm
might make iterated learning self-sustaining. Of course, it is easy to think of situations
where this feature would be highly desirable (eg, school teaching, social transmission of
norms, legends, jokes, etc.)
We show how keeping the length of the training sessions growing slightly
allows iterated learning to be sustained in perpetuity.

In the first part of the paper, we characterize iterated learnability in geometric terms
and show how a slight, steady increase in the lengths of the learning sessions
ensures self-sustainability for any discrete language class.
In the second part, we tackle the nondiscrete case and investigate
self-sustainability for iterated Bayesian linear regression.
In all cases, self-sustainability requires
making the underlying Markov process time-inhomogeneous
in order to stay out of equilibrium. 
This gives us an opportunity to offer a few thoughts on the growing importance
of non-equilibrium in natural algorithms.

\setlength{\columnsep}{17pt}%
\begin{wrapfigure}{r}{.48 \textwidth}
\vspace{-0.3cm}
\hspace{0.10cm}
\includegraphics[width=0.5 \textwidth, height=0.15 \textheight]{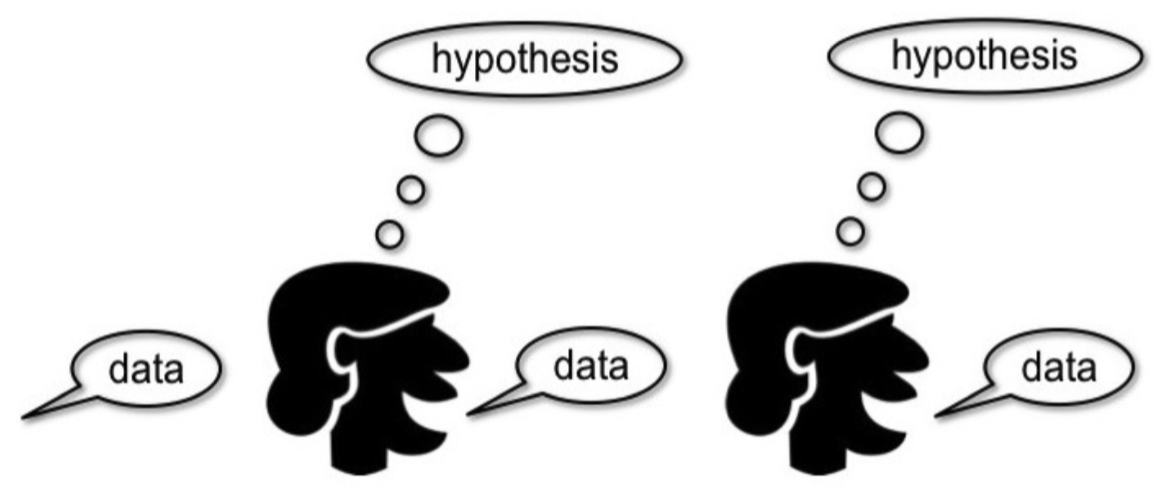} 
\par
\vspace{0.5cm}
\setlength\baselineskip{2.5ex}
{\footnotesize \hspace{1.4cm}  Fig.1: Chained iterated learning.}
\end{wrapfigure}

\paragraph{Background.}\
Following~\cite{beppu2009iterated,griffiths2007language,griffiths2005bayesian,rafferty2009convergence,kirby2014iterated,griffiths2008theoretical,perfors2011language,rafferty2014analyzing,smith2009iterated,kalish2007iterated},
we begin with {\em chained iterated learning}:
a learner's prior is modeled by a distribution over a hypothesis space ${\mathcal H}$, which is
itself equipped with a likelihood function:  ${\mathbb P}[d|h]$ indicates
the probability of generating data $d\in {\mathcal D}$ given 
the hypothesis $h\in {\mathcal H}$.
The first learner samples $m_1$ items {\em iid} from the 
initial hypothesis $\bm{h}_{\text{\tt init}}$: 
these items provide the training data ${\mathbf d}_1= (d_{1,1},\ldots, d_{1,m_1})$
with which the first learner Bayes-updates its prior. Its posterior is
given by setting $t=1$ in this formula:

\begin{equation}\label{bayesUpdate-t}
{\mathbb P}[h | {\mathbf d}_t] 
= {\mathbb P}[{\mathbf d}_t | h]\,  {\mathbb P}[h]/ {\mathbb P}[{\mathbf d}_t],
\hspace{0.3cm}
\text{with } 
{\mathbb P}[{\mathbf d}_t]=  
\sum_{h\in {\mathcal H}} {\mathbb P}[{\mathbf d}_t| h]  \, {\mathbb P}[h].
\end{equation}
From that point on, each successive learner updates its prior 
from their predecessor. 
For any $t>1$, learner $t$ receives $m_t$ items sampled
from the posterior of agent $t-1$ to form the training set ${\mathbf d}_{t}$.
To do that, she picks a random hypothesis $h$ from ${\mathcal H}$ with probability 
${\mathbb P}[h| {\mathbf d}_{t-1}] $ (the posterior of learner $t-1$)
and then samples $m_t$ items {\em iid} from $h$
to form ${\mathbf d}_t\in {\mathcal D}^{m_t}$.
The posterior ${\mathbb P}[h | {\mathbf d}_t]$ is derived according to~(\ref{bayesUpdate-t}).
Note that learner~$t$ has no direct access to the
posterior of learner $t-1$ but only to data drawn from a hypothesis
sampled from the posterior.
Our formulation assumes a discrete space ${\mathcal H}$ but extends to 
continuous settings, as we show in~\S\ref{sec:gaussian}.

In the case of linguistic transmission, 
each hypothesis $h\in {\mathcal H}$ 
is a ``knob" whose setting is given 
by a number between 0 and 1, specifically the prior probability ${\mathbb P}[h]$.
All learners share the same prior. 
Picking some $h$ from that prior specifies a {\em language} (also denoted $h$
for convenience). In this case, a language is defined as a probability
distribution over ${\mathcal D}$, interpreted here as a set of {\em sentences}.
In this way, the prior can be
viewed as a mixture over ${\mathcal H}$: by abuse of terminology,
we call it a {\em mixed} hypothesis, which we distinguish from 
a {\em pure} hypothesis of the form $h\in {\mathcal H}$ (corresponding to a single-point distribution). 
Access to language $h$ is achieved by random sampling:  
the sentence $d\in  {\mathcal D}$ is picked with 
probability ${\mathbb P}[d | h]$. 

Iterated learning proceeds as follows.
After selecting language $h$
with probability ${\mathbb P}[h|  {\mathbf d}_{t-1}]$,
learner $t$ collects $m_t$ independent samples from $h$.
Thus, given a tuple ${\mathbf d}_t=  (d_{1},\ldots, d_{m_t})$
of sentences from ${\mathcal D}$,  the likelihood ${\mathbb P}[ {\mathbf d_t} | h]$
is equal to $\prod_{1\leq k\leq m_t} {\mathbb P}[d_k | h]$.
The learner is now ready to Bayes-update its prior.
Of course, the first one $(t=1)$ samples directly from 
the language $\bm{h}_{\text{\tt init}}$ chosen for iterated learning.
The notation is boldfaced to indicate that $\bm{h}_{\text{\tt init}}$
may be a mixed hypothesis or, in other words, a distribution over hypotheses.

Suppose that ${\mathcal D}=\{d_1,\dots,d_s\}$ and ${\mathcal H}=\{h_1,\dots,h_n\}$ 
are both finite. 
While sampling from the posterior of learner $t-1$,
if learner $t$ winds up choosing~$h_i$ then, by Bayesian updating,
the probability $P^{|t}_{ij}$ that its posterior picks $h_j$ is given by:

\begin{equation}\label{Pij}
P^{|t}_{ij} 
=  
\sum_{{\mathbf d} \in {\mathcal D}^{m_t}} {\mathbb P}[h_j | {\mathbf d}]\,
{\mathbb P}[{\mathbf d} | h_i] 
=
\sum_{{\mathbf d} \in {\mathcal D}^{m_t}}
\frac{
{\mathbb P}[{\mathbf d} | h_i] \,  {\mathbb P}[{\mathbf d} | h_j]\,  {\mathbb P}[h_j]
}
{\sum_{k=1}^n {\mathbb P}[{\mathbf d} | h_k]  \, {\mathbb P}[h_k]
}\, .
\end{equation}
\smallskip

To our knowledge, the entire literature on the topic assumes
a common, fixed sample size for all the learners: $m_t=m$.
Equation~(\ref{Pij}) can be then interpreted as marginalizing a Gibbs sampler over the data space,
which creates a Markov chain over the hypothesis space ${\mathcal H}$:
if $\bm{h}^t$ denotes the row vector formed by the
$n$ probabilities ${\mathbb P}[h_k\,|\, {\mathbf d}_t]$, then
$\bm{h}^t= \bm{h}^{t-1}P^{t}$, where $\bm{h}^0= \bm{h}_{\text{\tt init}}$.
Assuming ergodicity (in this case, a fairly
inconsequential technical assumption), the chain can be shown to converge to a unique stationary distribution
$\bm{h}$. It can be easily checked that it coincides with the prior:
$\bm{h}= ({\mathbb P}[h_1], \ldots, {\mathbb P}[h_n])$~\cite{griffiths2005bayesian, norris1998markov};
see~\cite{rafferty2009convergence,rafferty2014analyzing} for an analysis of the mixing time
in specific linguistic scenarios.
This convergence reveals the long-term unsustainability of iterated learning. 
We show how diversifying the sample sizes $m_t$, hence making
the Markov chain time-inhomogeneous, can overcome this weakness.

\paragraph{Our results.}\
In~\S\ref{sec:discrete}, we show how to achieve self-sustainability
in the discrete setting~\cite{griffiths2007language,griffiths2005bayesian}, using
only a logarithmically increasing sample size; specifically, the new hypothesis to be learned
is acquired by all the (infinitely many) learners with probability 
at least $1-\eps$ using a sample size of $O( \log \frac{t}{\eps})$ for
the $t$-th learner. The constant factor depends on the geometry of
the hypothesis space.  By relaxing the objective and allowing learners to
settle on an arbitrarily close approximation of the hypothesis to be learned,
we can remove all dependency on the geometry of the hypothesis space.

In~\S\ref{sec:gaussian}, we extend the iterated learning model
to a Gaussian setting for an infinite hypothesis space and show that
a sample size of $O(t)^{1+ o(1)}$  is sufficient to ensure self-sustainability.
We also show that allowing learners to pick their teachers at random cuts down
the sample size to $O(\log t)^{1+o(1)}$.  The arguments used for the discrete case
bump into singularities so we use a different approach, which allows us to exploit various
``stability" properties of the Gaussian setting.

In~\S\ref{sec:linear}, we turn our attention to the iterated version of Bayesian linear regression
and prove a high-probability statement about self-sustainability. This requires spectral arguments
from random matrix theory and, in particular, bounds on the lowest singular
value of Wishart matrices. 

\paragraph{Discussion.}\
Before moving to the technical part of this work, we add a few thoughts about
its larger context and relevance.  
For a dynamicist, the loss of Quenya 
is a byproduct of the memory-erasing ergodicity implied by mixing.
For a physicist,  the loss is due to the Second Law of thermodynamics
and the bounded supply of free energy available
to each agent: together these two constraints make it impossible to keep the system out of equilibrium.
For a biologist, this entropic pull toward equilibrium is the hallmark of a dying system.  
Evolution is nature's
attempt to optimize the absorption of free energy into work while
maximizing the production of entropy. The first requirement is keeping
the system out of equilibrium over timescales well in excess of the metabolic rate
(here, the teaching rate).  From that perspective,
our work can be seen as an effort to find out the minimum conditions
necessary to keep a target dynamics active in perpetuity. 
There are several approaches to this question
and the two we follow are among the simplest:
(i) increasing the supply of free energy (eg, lengthening the 
training sessions) and (ii) mixing timescales (eg, rewiring the communication network).

Most of the work on Markov chains in theoretical computer science
regards mixing as a blessing: large spectral gaps are good while small ones are
to be avoided.  In biology, however, mixing often means death. In fact, much of life can be seen
as nature's attempt to keep mixing at bay.  This paper explores what can be
done to prevent a Markov chain from reaching equilibrium. We expect this theme
to gain prominence in future work on natural algorithms.

\section{Self-Sustainability} \label{sec:discrete}

We show how to make iterated learning self-sustaining in the presence
of a finite hypothesis space ${\mathcal H}=\{h_1,\dots,h_n\}$. 
This involves specifying a sequence of training session lengths $m_1,m_2,\ldots$ so that
the posterior of any learner ends up differing from  $\bm{h}_{\text{\tt init}}$ 
by an arbitrarily small amount.   Formally,
given any $\delta, \eps\geq 0$, we say that iterated learning is {\em $(\delta,\eps)$-self-sustaining}
if, for any learner,  with probability at least $1-\eps$, a random $h\in {\mathcal H}$ picked
from its posterior distribution differs from $\bm{h}_{\text{\tt init}}$ in total variation
by at most $\delta$. 
We recall a few facts:
the hypothesis $h$ denotes a language modeled as a probability distribution over ${\mathcal D}$;
the total variation distance is half the $\ell_1$-norm; and 
the posterior of learner $t$ after the $t$-th iteration
is defined by marginalizing ${\mathbb P}[h | {\mathbf d}_t]$ over all samples ${\mathbf d}_t$ 
drawn from a random $h$ picked from the posterior of learner $t-1$
(or $\bm{h}_{\text{\tt init}}$ if $t=1$).
As a shorthand, we speak of $\eps$-self-sustainability to refer to the case $\delta=0$.

The parameters $\delta$ and $\eps$ allow us to distinguish between
two metrics: the distance between two languages over ${\mathcal D}$ and the distance
between two mixtures over  ${\mathcal H}$.
The two notions could differ widely. For example, if all of ${\mathcal H}$ corresponds to
languages very close to $\bm{h}_{\text{\tt init}}$, to achieve
$(\delta,\eps)$-self-sustainability might be easy for a tiny $\delta>0$ 
but hopelessly difficult for $\delta=0$.
The complexity of iterated learning depends on 
the geometry of the languages formed by the pure hypotheses.
This is best captured by introducing a metric that, though more specialized
than the total variation (it works only on the simplex of
probability vectors) brings all sorts of technical benefits:
the {\em root-sine distance} between two probability distributions
$\bm{a}=(a_1,\ldots, a_s)$ 
and $\bm{b}= (b_1,\ldots, b_s)$ over ${\mathcal D}$ is defined as

\begin{equation}\label{d_RS}
d_{RS}(\bm{a},\bm{b})
= \, \sqrt{ \frac{1}{2} \sum_{i,j=1}^s
\left (\sqrt{ a_ib_j}- \sqrt{a_jb_i}\,
\right)^2 }
=  \sqrt{ 1 -  \Bigl( \sum_{i=1}^s \sqrt{ a_ib_i}\Bigr)^2}\, .
\end{equation} 

\medskip

It would be surprising if this distance had not been used before, but we 
could not find a reference.  
We prove that it is indeed a metric in the Appendix and also explain its name.
We show that it is related to the
Hellinger, Bhattacharyya
and total variation distances, 
$d_H$, $d_B$, $d_{TV}$ by the following relations:

\begin{equation}\label{DistCompare}
\begin{cases}
\, d_H= \sqrt{1-\sqrt{1-d_{RS}^2}}\, ; \vspace{.2cm} \\
\, d_B= -\frac{1}{2}\ln (1- d_{RS}^2)\, ; \vspace{.2cm} \\
\, d_{TV}\leq \sqrt{2s} \, d_{RS}.
\end{cases}
\end{equation}

\subsection{The results}

We focus on the ``pure" case $\bm{h}_{\text{\tt init}}\in {\mathcal H}$,
and later briefly discuss how to generalize the method to mixed hypotheses.
Using the shorthand ${\tt d}_{ij}$ for $d_{RS}({\mathbb P}[\cdot |h_i],{\mathbb P}[\cdot |h_j])$,
we define ${\tt d}_i: =\min_{j: j\neq i} {\tt d}_{ij}$.
Let $\bm{p}=(p_1,p_2,\ldots, p_n)$ be the prior distribution over ${\mathcal H}$,
where $p_i:={\mathbb P}[h_i]$.
We can obviously assume that each $p_i$ is positive and
that all the pure hypotheses are distinct, hence ${\tt d}_i>0$.
The two theorems below assume that $\bm{h}_{\text{\tt init}}=h_1$.

\medskip

\begin{theorem}\label{discreteEps}
$\!\!\! .\,\,$
For any positive $\eps <1$, the following sample size sequence makes
iterated learning $\eps$-self-sustaining:
$$
m_t= \frac{4}{{\tt d}_1^2} 
\ln \frac{nt}{ \eps \, p_1} 
= 
\frac{4}{{\tt d}_1^2} \Bigl( \log \frac{t}{\eps}  + C \Bigr),
$$
for some $C>0$ independent of $t, \eps, {\tt d}_1$.
\end{theorem}

\bigskip
\noindent
The factor~$4$ can be reduced to $2^{1+o(1)}$ if we adjust the constant $C$.
It is to be expected that the lengths of the training sessions 
should grow to infinity as $p_1$ tends to zero,
as the vanishing prior makes it increasingly difficult for the posteriors to ``attach" to $h_1$.
The session lengths are sensitive to the minimum distance 
between the languages specified by ${\mathcal H}$ and the target language $h_1$.
Settling for $(\delta, \eps)$-self-sustainability allows us to remove this dependency.
\bigskip

\begin{theorem}\label{discreteDeltaEps}
$\!\!\! .\,\,$
For any positive $\delta,  \eps <1$,
the following sample size sequence makes
iterated learning $(\delta,\eps)$-self-sustaining:
$$
m_t= \frac{8s n^2}{\delta^2} 
\Bigl( \ln \frac{t}{ \eps } +C\Bigr).
$$
for some $C>0$ independent of $t, \delta, \eps$. 
\end{theorem}

\bigskip

\subsection{The proofs}

To establish Theorem~\ref{discreteEps},
we begin by bounding the off-diagonal elements of the transition matrix $P_{ij}^{|t}$
for $i\neq j$.  By~(\ref{Pij}) and Young's inequality,
\begin{equation*} 
\begin{split}
P^{|t}_{ij} 
& \leq  \sum_{{\mathbf d} \in {\mathcal D}^{m_t}}
\frac{
{\mathbb P}[{\mathbf d} | h_i] \,  {\mathbb P}[{\mathbf d} | h_j] p_j
}
{{\mathbb P}[{\mathbf d} | h_i] p_i+  {\mathbb P}[{\mathbf d} | h_j] p_j
}
=    \frac{p_j}{p_i} \sum_{{\mathbf d} \in {\mathcal D}^{m_t}}
\frac{
\bigl(\frac{p_i}{p_j}\bigr){\mathbb P}[{\mathbf d} | h_i] \,  {\mathbb P}[{\mathbf d} | h_j]
}
{\bigl(\frac{p_i}{p_j}\bigr){\mathbb P}[{\mathbf d} | h_i] +  {\mathbb P}[{\mathbf d} | h_j]
} \\
&\leq
\frac{1}{2}\sqrt{\frac{p_j}{p_i}}
\sum_{{\mathbf d} \in {\mathcal D}^{m_t}}
\sqrt{
{\mathbb P}[{\mathbf d} | h_i] \,  {\mathbb P}[{\mathbf d} | h_j] }
= \frac{1}{2}\sqrt{\frac{p_j}{p_i}}
\Bigl (\sum_{d \in {\mathcal D} }
\sqrt {
{\mathbb P}[d | h_i]\,
{\mathbb P}[d | h_j]  } \,
\Bigr)^{m_t}\\
&\leq 
\frac{1}{2}\sqrt{\frac{p_j}{p_i}}
\exp{
\biggl\{ \frac{m_t }{2}
\Bigl( \Bigl( 
\sum_{d \in {\mathcal D} }
\sqrt {
{\mathbb P}[d | h_i]\,
{\mathbb P}[d | h_j]  } \,
\Bigr)^2  -1 \Bigr) \biggr\} }
 \, .
\end{split}
\end{equation*}
By definition of the root-sine distance, we have

\begin{equation}\label{Pij2}
P^{|t}_{ij} 
\leq \,
\frac{1}{2}\sqrt{\frac{p_j}{p_i}}
e^{-\frac{1}{2}{\tt d}_{ij}^2 m_t }\ \ \ \ 
(i\neq j).
\end{equation}
\medskip
Setting $i=1$ in~(\ref{Pij2}) and summing over $2\le j\le n$, it follows by
Cauchy-Schwarz that

\begin{equation}\label{t-ij<}
\sum_{j=2}^n  P^{|t}_{1j} 
\leq \, \frac{1}{2}
\sqrt{\frac{n(1-p_1)}{p_1}} \,
e^{-\frac{1}{2}{\tt d}_{1}^2 \, m_t } \, .
\end{equation}

\begin{lemma}\label{Pij-UB}
$\!\!\! .\,\,$
Let $ P^{\leq t}$ denote the matrix product $P^{|1}\cdots P^{|t}$.
For any $t\geq 1$,
$$ 
\sum_{j=2}^n  P^{\leq t}_{1j}
\leq \, \frac{1}{2}
\sqrt{\frac{n(1-p_1)}{p_1}} \, \sum_{s=1}^t  e^{-\frac{1}{2}{\tt d}_{1}^2 \, m_s }\, .$$
\end{lemma}
\proof
The case $t=1$ follows from~(\ref{t-ij<}), so we assume that $t>1$ 
and proceed by induction.
By~(\ref{t-ij<}),
\begin{equation*}
\begin{split}
\sum_{j=2}^n  P^{\leq t}_{1j}
&= \sum_{j=2}^n  \sum_{k=1}^n P^{\le t-1}_{1k}\, P^{|t}_{kj} 
=   P^{\leq t-1}_{11} \sum_{j=2}^n  P^{|t}_{1j}
+   \sum_{j,k=2}^n P^{\le t-1}_{1k} \, P^{|t}_{kj} \\
&\leq 
 \sum_{j=2}^n  P^{|t}_{1j}
+   \sum_{k=2}^n P^{\le t-1}_{1k} 
\leq
\sqrt{\frac{n(1-p_1)}{4p_1}} \,
e^{-\frac{1}{2}{\tt d}_{1}^2 \, m_t }
+  \sqrt{\frac{n(1-p_1)}{4p_1}} \, \sum_{s=1}^{t-1} e^{-\frac{1}{2} {\tt d}_{1}^2 \, m_s }.
\end{split}
\end{equation*}
\hfill $\Box$
\proofend

\noindent
Given $0< \eps < 1$, we constrain the sequence $(m_t)$ to satisfy:
\begin{equation}\label{cond-nt}
 \sum_{s=1}^t  e^{-\frac{1}{2}{\tt d}_{1}^2 m_s } <  \eps \sqrt{\frac{4p_1}{n(1-p_1)}}\, .
\end{equation}
For example, we can pick the sequence
\begin{equation*}
m_t=  \frac{1}{{\tt d}_{1}^2}\ln \frac{n(1-p_1)t^4}{ \eps ^2p_1}\, .
\end{equation*}
A closer look at the calculation shows
that the factor $t^4$ can be reduced to $C_\alpha t^{2+\alpha}$ for any small $\alpha>0$
and a suitable constant $C_\alpha>0$, which makes the dependency on $t$ arbitrarily close 
to $(2/ {\tt d}_{1}^2)\ln t$.
After the $t$-th iteration, the posterior of the $t$-learner becomes
$\bm{h}^t= \bm{h}^0  P^{\leq t}$, where $\bm{h}^0 = (1,0,\ldots, 0)$.
By Lemma~\ref{Pij-UB} and~(\ref{cond-nt}),
$$ 
{\mathbb P}[h= h_1| {\mathbf d}_t ]=
\bm{h}^t_1 = P^{\leq t}_{11}= 
1- \sum_{j=2}^n  P^{\leq t}_{1j}
\geq  1- 
 \sqrt{\frac{n(1-p_1)}{4p_1}} \sum_{s=1}^t  e^{-{\tt d}_{1}^2 m_s } > 1- \eps .
$$
\hfill $\Box$
\proofend

To prove Theorem~\ref{discreteDeltaEps},
we set a target distance $\rho:= \delta/(n\sqrt{2s})$
and find a subset $A\subseteq {\mathcal H}$
such that 
(i) ${\tt d}_{1j}\leq \rho n$ for $j\in A$ and
(ii) ${\tt d}_{ij}\geq \rho$ for $i\in A$ and $j\not\in A$.
To see why such a subset must exist, consider spheres centered 
at $\bm{h}_{\text{\tt init}}=h_1$ of radius $k\rho$, for $k=1,\ldots, n+1$
(with respect to $d_{RS}$). These define $n+1$ disjoint (open) regions
and, by the pigeonhole principle, at least 
one of them must be empty. We set $A$ to include 
all the points in the regions 
preceding the empty one; note that $h_1\in A$. The claim follows from
the triangular inequality. 
We begin with a straightforward generalization of Lemma~\ref{Pij-UB}:

\begin{equation}\label{PAB-UB}
\sum_{j\not\in A}  P^{\leq t}_{1j}
\leq \, \frac{1}{2}
\sqrt{\frac{n(1-p_A)}{p_A}}
            \, \sum_{s=1}^t  e^{-\frac{1}{2}\rho ^2 \, m_s }\,  ,
\end{equation}
where $p_A= \min_{i\in A}p_i$
To see why, first observe that, by~(\ref{Pij2}) and Cauchy-Schwarz,
for any $i\in A$,

\begin{equation}\label{t-AB}
\sum_{j\not\in A}  P^{|t}_{ij} 
\leq \, \frac{1}{2}
\sqrt{\frac{n(1-p_i)}{p_i}} \,
e^{-\frac{1}{2}\rho^2 \, m_t } \, .
\end{equation}
Inequality~(\ref{PAB-UB}) is proven by induction:
the case $t=1$ follows from~(\ref{t-AB}); for $t>1$,

\begin{equation*}
\begin{split}
\sum_{j\not\in A}  P^{\leq t}_{1j}
&= \sum_{j \not\in A}  \sum_{k=1} P^{\leq t-1}_{1k}\, P^{|t}_{kj} 
\leq
\max_{i\in A}  \sum_{j\not\in A}  P^{|t}_{ij} 
+ 
\sum_{j\not\in A}  P^{\leq t-1}_{1j} \\
&\leq
\max_{i \in A}
\sqrt{\frac{n(1-p_i)}{4p_i}} \,
e^{-\frac{1}{2}\rho^2 \, m_t }
+  \sqrt{\frac{n(1-p_A)}{4p_A}} \, 
           \sum_{s=1}^{t-1} e^{-\frac{1}{2} \rho^2 \, m_s }\, ,
\end{split}
\end{equation*}
which proves~(\ref{PAB-UB}).
Recall that $\bm{h}^t= \bm{h}^{t-1}P^{t}$,
where $\bm{h}^0=  (1,0,\ldots, 0)$.
Given a random $h$ from ${\mathcal H}$,
$$
{\mathbb P}[h\in A | {\mathbf d}_t ]
= 
\sum_{i\in A} \bm{h}^t_i 
=   \sum_{j\in A}  P^{\leq t}_{1j}
=  1- \sum_{j\not\in A}  P^{\leq t}_{1j}
\geq
1- \sqrt{\frac{n(1-p_A)}{4p_A}} 
             \sum_{s=1}^t  e^{-\frac{1}{2}\rho^2 m_s } \, .
$$
Setting
\begin{equation}\label{defn-ntAB}
m_t=  \frac{1}{\rho^2}\ln \frac{n(1-p_A)t^4}{ \eps ^2p_A}\, 
\end{equation}
ensures that ${\mathbb P}[h\in A | {\mathbf d}_t ]> 1- \eps$.
The root-sine distance between the languages denoted by $h_1$ and any $h\in A$
is at most $\rho n$, so that, by~(\ref{DistCompare}), the total variation distance is bounded
by $\sqrt{2s} \rho n= \delta$, which concludes 
the proof of Theorem~\ref{discreteDeltaEps}.
\hfill $\Box$
\proofend

So far, we have analyzed only the ``pure" case $\bm{h}_{\text{\tt init}}\in {\mathcal H}$.
The idea of the training is to prevent the prior to ``drag" the posterior mixture all
across ${\mathcal H}$.  It should be clear that a similar result obtains if
$\bm{h}_{\text{\tt init}}\in \Delta {\mathcal H}$ is concentrated on a subset $A$ of ${\mathcal H}$.
The proof follows the path charted in Theorem~\ref{discreteDeltaEps}
and need not be repeated here.  It is crucial to note, however, that this result is to be understood
in a coarse-graining sense:  iterated learning cannot ensure that the original weights in the mixture
$\bm{h}_{\text{\tt init}}$ are retained but only that $A$ contributes most of the mass
in the posteriors. To retain the weights would require changing the
stationary distribution to conform with $\bm{h}_{\text{\tt init}}$, as the process
unfolds, something that straightforward Bayesian learning seems unable to do.
Learning pure hypotheses bypasses that difficulty.

\subsection{Applications}

We briefly discuss a direct application of our results to a well-known model of language acquisition
via iterated learning and we mention some natural extensions of the techniques.

\paragraph{Language evolution.}\ 
Rafferty et al.~\cite{rafferty2009convergence} show how iterated learning fails rapidly
in a simple model of language evolution. Given $n$ hypotheses, iterated learning with
fixed-length training sessions ceases to learn anything new after only $O(\log n\log\log n)$ rounds.
The previous theorems show how to turn this around and achieve self-sustainality.
In the model, ${\mathcal H}= \{h_1,\ldots, h_n\}$, where $n=2^k$ and $h_i$ denotes the language
whose sentences are words in $\{0,1,?\}^k$ with exactly $m$ question marks and $0,1$ matching
the binary decomposition of $i-1$ outside the question marks.
For example, if $k=4$ and $m=2$, then $h_3$ denotes the language
$$\{\,  00??, 0?1?, ?01?, 0??0, ?0?0, ??10 \, \}.$$
We can assume that $m$ is much smaller than $k$. 
Each language has the same length $\binom{k}{m}$ and the total
number of sentences is $s= \binom{k}{m} 2^{k-m}$. The prior is given by
${\mathbb P}[h_i] = p_i= 1/n$. Given a hypothesis $h_i$, 
${\mathbb P}[d | h_i] = 1/\binom{k}{m}$ if $d$ has $m$ question marks
and match the bits of $i-1$ elsewhere; else it is 0 (and $d, h$ are called incompatible).
Given $h\in {\mathcal H}$,

\begin{equation*}
\begin{cases}
\,    {\mathbb P}[d] = \sum_{h\in {\mathcal H}} {\mathbb P}[d | h] {\mathbb P}[h]= 2^{m-k}/ \binom{k}{m}
     \, ; \vspace{.2cm}   \\
\,    {\mathbb P}[h|d]=  {\mathbb P}[d|h] {\mathbb P}[h] / {\mathbb P}[d] = 2^{-m}
\hspace{.5cm}
\text{(or $0$ if $d, h$ are incompatible)}.
\end{cases}
\end{equation*}
We easily check that 
${\tt d}_1^2
=    1 -  \bigl( \sum_{i=1}^s \sqrt{ a_ib_i}\, \bigr)^2
\geq   1 -  \bigl(\frac{m}{k}\bigr)^2 > \frac{1}{2}$; hence,
by Theorem~\ref{discreteEps}, session lengths $m_t$ no larger than
$O( \log \frac{t}{\eps})$ are sufficient to maintain $\eps$-self-sustainability.

\paragraph{Meanings and utterances.}\
In the use of iterated learning for
studying language evolution~\cite{griffiths2005bayesian,perfors2011language},
it is common to model the data ${\mathbf d}$ as a joint distribution
$(\bm{x},\bm{y})$ over a product space ${\mathcal X}^{m_t}\times {\mathcal Y}^{m_t}$.
The idea is to distinguish between ``meanings" $\bm{x}$ and ``utterances" $\bm{y}$.
In this setting,
${\mathbb P}[{\mathbf d} | h] = {\mathbb P}[\bm{y} | \bm{x}, h] \mu(\bm{x})$, 
where $\mu(\bm{x})$ is the probability of generating $\bm{x}$.
The transition matrix of the Markov chain thus becomes

\begin{equation}\label{PijNew}
\begin{split}
P^{|t}_{ij} 
&=  
\sum_{\bm{x} \in {\mathcal X}^{m_t}} \sum_{\bm{y} \in {\mathcal Y}^{m_t}} {\mathbb P}[h_j | \bm{x},\bm{y}]\,
{\mathbb P}[\bm{y} | \bm{x}, h_i] \mu(\bm{x}) \\
&=
\sum_{\bm{x} \in {\mathcal X}^{m_t}} \sum_{\bm{y} \in {\mathcal Y}^{m_t}}
\frac{
{\mathbb P}[\bm{y} | \bm{x}, h_i] \,  {\mathbb P}[\bm{y} |\bm{x}, h_j]\,  {\mathbb P}[h_j]
}
{\sum_{k=1}^m  {\mathbb P}[\bm{y} |\bm{x}, h_k]\,  {\mathbb P}[h_k] }
\mu(\bm{x})\, .
\end{split}
\end{equation}

\medskip

\noindent
Since the output $\bm{y}$ now depends on both the hypothesis and the input data, 
we redefine ${\tt d}_{ij}$ as the root-sine distance between the two distributions
${\mathbb P}[\bm{y}|\bm{x},h_i]\mu(\bm{x})$ and ${\mathbb P}[\bm{y}|\bm{x},h_j]\mu(\bm{x})$:
\begin{equation}\label{gap_2}
{\tt d}_{ij}' :=1-\left(\sum_{\bm{x}\in{\mathcal X}}\sum_{\bm{y}\in{\mathcal Y}}\sqrt{{\mathbb P}[\bm{y}|\bm{x},h_i]{\mathbb P}[\bm{y}|\bm{x},h_j]}\, \mu(\bm{x})\right)^2
\end{equation}
and we define ${\tt d}_i': =\min_{j: j\neq i} {\tt d}_{ij}'$.
Given any $i\neq j$,

\begin{equation*}
\begin{split}
P^{|t}_{ij} 
& 
\leq
\sum_{\bm{x} \in {\mathcal X}^{m_t}} \sum_{\bm{y} \in {\mathcal Y}^{m_t}}
\frac{
{\mathbb P}[\bm{y} | \bm{x}, h_i] \,  {\mathbb P}[\bm{y} |\bm{x}, h_j]\,  p_j
}
{{\mathbb P}[\bm{y} |\bm{x}, h_i]\,  p_i+{\mathbb P}[\bm{y} |\bm{x}, h_j]\,  p_j}
\mu(\bm{x})
\\
&\leq
\frac{1}{2}\sqrt{\frac{p_j}{p_i}}
\sum_{\bm{x} \in {\mathcal X}^{m_t}} \sum_{\bm{y} \in {\mathcal Y}^{m_t}}
\sqrt{
{\mathbb P}[\bm{y} | \bm{x}, h_i] \,  {\mathbb P}[{\mathbf y} | \bm{x},h_j] }\, \mu(\bm{x}) \\
&\leq \frac{1}{2}\sqrt{\frac{p_j}{p_i}}
\Bigl (\sum_{\bm{x} \in {\mathcal X}} \sum_{\bm{y} \in {\mathcal Y}}
\sqrt {
{\mathbb P}[\bm{y} | h_i]\,
{\mathbb P}[\bm{y} | h_j] }\,  \mu(\bm{x}) \,
\Bigr)^{m_t}\\
&\leq 
\frac{1}{2}\sqrt{\frac{p_j}{p_i}}
\exp{
\biggl\{ \frac{m_t}{2} 
\Bigl( \Bigl( 
\sum_{\bm{x} \in {\mathcal X}} \sum_{\bm{y} \in {\mathcal Y}}
\sqrt {
{\mathbb P}[\bm{y} |\bm{x}, h_i]\,
{\mathbb P}[\bm{y} |\bm{x}, h_j]} \,\mu(\bm{x}) \,
\Bigr)^2  -1 \Bigr) \biggr\} }
 \, .
\end{split}
\end{equation*}
This gives us this new version of inequality~(\ref{Pij2}),
which we can use as the basis for a repeat of the argument 
of the previous section:

\begin{equation}
P^{|t}_{ij} 
\leq \,
\frac{1}{2}\sqrt{\frac{p_j}{p_i}}
e^{-\frac{1}{2}{\tt d}_{ij}'^2 m_t } \ \ \ \ 
(i\neq j).
\end{equation}

\section{Iterated Learning in Continuous Spaces}\label{sec:gaussian}

When iterated learning operates over a hypothesis space ${\mathcal H}$ parametrized
continuously, say, in ${\mathbb R}$,  
the minimum root-sine distance usually vanishes and the previous
arguments run into singularities and collapse.  A new approach is needed.
To make our discussion concrete,
we assume that the prior distribution of each learner is a 
Gaussian ${\mathbb P}[h]\sim N(\bar{\mu},\bar{\sigma}^2)$
and that the likelihood of producing data $d$ given hypothesis $h$ 
is also normal:  ${\mathbb P}[d|h]= N(h,\sigma^2)$. 
The likelihood can also be understood as a noisy measurement of $h$: $d=h+\phi$, where the noise $\phi\sim N(0,\sigma^2)$.
We assume that the data received by the first learner 
comes from $N(\mu_0,\sigma_0^2)$.  This is the simplest instance of a continuous
setting in which the root-sine distance argument fails. We discuss it in some detail, considering
both chained learning and its generalizations;
and then we use the results to treat the case of iterated Bayesian linear regression.

During its training session, the $t$-th learner receives 
data ${\mathbf d}_t=(d_{t,1},\dots,d_{t,m_t})$ from its predecessor: it is obtained
by first picking a random hypothesis $h$ from the posterior of learner $t-1$
and then collecting $m_t$ independent random samples from $N(h,\sigma^2)$.
For the case $t=1$, we can treat the original teacher as learner $0$ with its posterior 
equal to $N(\mu_0,\sigma_0^2)$.
Learner $t$ Bayes-updates its posterior as follows: 
$$
{\mathbb P}[h|{\mathbf d}_t]\propto {\mathbb P}[{\mathbf d}_t|h]{\mathbb P}[h]\propto
\exp\Bigl(-\frac{1}{2\sigma^2}\sum_{i=1}^{m_t}{(d_{t,i}-h)^2}\Bigr)\exp\Bigl(-\frac{1}{2\bar{\sigma}^2}(h-\bar{\mu})^2\Bigr),$$
which is still Gaussian, with mean and variance denoted by
$\mu_t$ and $\sigma_t^2$, respectively. 
Carrying out the usual square completion gives up these update rules: for $t>0$,
\begin{equation}\label{normal_update}
\begin{cases}
\,\,  \mu_t = \frac{1}{\bar{\tau}+m_t\tau}  \left(\bar{\tau}\bar{\mu}+\tau(d_{t,1}+d_{t,2}+\dots+d_{t,m_t})\right)
       \vspace{.2cm} \\
\,\,
\tau_t \, = \bar{\tau}+m_t\tau ,
\end{cases}
\end{equation}
where we define the precisions $\tau=1/\sigma^2$, $\bar{\tau}=1/\bar{\sigma}^2$, and $\tau_t=1/\sigma_t^2$.
We say that iterated learning is {\em $\eps$-self-sustaining}
if $|{\mathbb E}\, \mu_t - \mu_0|\leq \eps$ and
$\sigma_t^2 + \Varr \mu_t $ remains bounded for all $t$.
If $\sigma_t^2 + \Varr \mu_t  \rightarrow 0$ as $t\rightarrow \infty$,
we say that iterated learning is {\em strongly $\eps$-self-sustaining}.
We consider successively the case of chained iterated learning and the more challenging 
``hopping" scenario 
in which a new learner picks a random teacher from the past (instead of the previous one).

\subsection{Chained learning }

In chained iterated
learning, the data $d_{t,i}$ is a noisy message drawn 
from the posterior of the $(t-1)$-th learner; hence $d_{t,i}\sim N(\mu_{t-1},\sigma_{t-1}^2+\sigma^2)$.
In view of \eqref{normal_update}, $\mu_t$ is itself Gaussian.
By taking the expectation and variance of equation \eqref{normal_update}, we find
the following recursive relations for 
${\mathbb E}\, \mu_t $ and $\Varr \mu_t$: for $t>0$,
\begin{equation}\label{normal_update2}
\begin{cases}
\,\,    {\mathbb E}\, \mu_t = \frac{1}{\bar{\tau}+m_t\tau}
\bigl(\bar{\tau}\bar{\mu}+m_t\tau \, {\mathbb E}\, \mu_{t-1}\bigr); \vspace{.3cm} \\
\,\, \Varr  \mu_t = \frac{m_t\tau^2}{(\bar{\tau}+m_t\tau)^2} \bigl(\Varr\, \mu_{t-1}
+\sigma_{t-1}^2+\sigma^2\bigr).
\end{cases}
\end{equation}

\medskip
\noindent
If we define $\beta_t:=m_t\tau/(\bar{\tau} + m_t\tau)$,
then \eqref{normal_update2} becomes ${\mathbb E}\, \mu_t
=\beta_t \, {\mathbb E}\, \mu_{t-1} +(1-\beta_t)\bar{\mu}$.
If $m_t=m$ is a constant, then so is $\beta_t$, and the recursive relation \eqref{normal_update2} becomes
$${\mathbb E}\, \mu_t -\bar{\mu}=\beta_1^t(\mu_0-\bar{\mu}),$$
which shows that ${\mathbb E}\, \mu_t$ converges to $\bar{\mu}$ exponentially fast.
As in the discrete case, iterated learning is not self-sustainable
with constant-length training sessions.
By letting $m_t$ increase as $O(t^{1+o(1)})$ order, however, we can achieve self-sustainability:

\begin{theorem}\label{th:chained}
For  any $0< \eps <1$, the following sample size sequence makes
chained iterated learning strongly $\eps$-self-sustaining:
$$m_t=  \frac{|\mu_0-\bar{\mu}|}
     { \eps }\Bigl(1+\frac{1}{c} \Bigr) \Bigl(\frac{\sigma}{\bar\sigma}\Bigr)^{\!2}
         \, t^{1+c},
$$
for an arbitrarily small constant $c>0$. 
\end{theorem}

\proof
We observe that
${\mathbb E}\, \mu_t$ is a convex combination of $\bar{\mu}$ 
and ${\mathbb E}\, \mu_s$ $(s<t$); specifically,
\begin{equation}
{\mathbb E}\, \mu_t=\prod_{s=1}^t\beta_s\mu_0+ \biggl(1-\prod_{s=1}^t\beta_s\biggr)\bar{\mu}.
\end{equation}
Because $\sum_{s>0}(1/s)^{1+c}< 1+ \int_1^\infty x^{-1-c}\, dx = 1+1/c$, we have
\begin{eqnarray*}
1&\geq & \prod_{s=1}^t\beta_s = \prod_{s=1}^t \left(1-\frac{\bar{\tau}}{m_s\tau+\bar{\tau}} \right)\ge 1-\sum_{s=1}^{t} \frac{\bar{\tau}}{m_s\tau+\bar{\tau}}
\\&\ge& 1- 
\frac{ \eps }{|\mu_0-\bar{\mu}|} \Bigl( \frac{c}{c+1} \Bigr) \sum_{s=1}^{\infty}\frac{1}{s^{1+c}} 
> 1-\frac{ \eps }{|\mu_0-\bar{\mu}|} \, .
\end{eqnarray*}
This shows that
$$
\left| {\mathbb E}\, \mu_t-\mu_0\right|
= \Bigl( 1-\prod_{s=1}^t\beta_s \Bigr) |\bar{\mu}-\mu_0|
\le  \eps .
$$
By~\eqref{normal_update}, $\sigma^2_t=1/\tau_t<1/m_t\tau \rightarrow 0$.
Since $\sigma_{t-1}^2\le \bar{\sigma}^2$ for $t>1$, it follows from~\eqref{normal_update2} that
$\Varr  \mu_t \le(\Varr  \mu_{t-1} +\sigma^2+\bar{\sigma}^2)/m_t$ for $t>1$,
and $\Varr  \mu_1 \le(\sigma_0^2+ \sigma^2)/m_1$.
Writing $M_t:= m_tm_{t-1}\dots m_1$, we have
\begin{equation*}
\begin{split}
M_t\Varr  \mu_t  
&\leq  M_{t-1}\Varr  \mu_{t-1} +M_{t-1}(\sigma^2 +\bar{\sigma}^2)\\
&\leq tM_{t-1}(\sigma_0^2+ \sigma^2+ \bar{\sigma}^2),
\end{split}
\end{equation*}
and thus $\Varr  \mu_t \le (\sigma_0^2+ \sigma^2+ \bar{\sigma}^2)t/m_t\rightarrow 0$
since $m_t =\Omega(t^{1+c})$.
\hfill $\Box$
\proofend

\subsection{Hopped learning}

We consider the ``hopped learning" scenario in which learner $t$ hops back to
pick a teacher from $\{0,1,\dots, t-1\}$ at random, and then samples $m_t$ bits of data from her
posterior.
The recursive relation for $\mu_t$ becomes
\begin{equation}\label{random_update}
\mu_t=\frac{\beta_t}{m_t}\sum_{s=0}^{t-1}\chi_{t,s}\sum_{i=1}^{m_t}d_{t,s,i}+(1-\beta_t)\bar\mu,
\end{equation}
where, given $t$, the random variable $\chi_{t,s}$ is 1 for a value of $s$ picked
at random between $0$ and $s-1$, and is zero elsewhere;
recall that $\beta_t:=m_t\tau/(\bar{\tau} + m_t\tau)$.
Hopped iterated learning provides access to earlier data, so one would expect the
lengths of the training sessions to grow more slowly than in chained learning.
The change is indeed quite dramatic:

\begin{theorem}
For any positive $\eps < |\mu_0-\bar{\mu}|$, the following sample size sequence makes
hopped iterating learning $\eps$-self-sustaining:
$$
m_t= B_c \frac{|\mu_0-\bar{\mu}|}{ \eps }\Bigl(\frac{\sigma}{\bar\sigma}\Bigr)^{\!2} (1+ \log t)^{1+c},$$ 
for an arbitrarily small $c>0$ and a constant $B_c$ that depends only on $c$.
\end{theorem}
\proof
By taking expectation on both sides of \eqref{random_update}, for any $t>0$,
$${\mathbb E}\, \mu_t=\frac{\beta_t}{t}\sum_{s=0}^{t-1}{\mathbb E}\, \mu_s+(1-\beta_t)\bar\mu,$$
We define $\gamma_1= \beta_1$ and, for $t>1$,
$$\gamma_t:=\left(1+\beta_1\right)\left(1+\frac{\beta_2}{2}\right)\cdots\left(1+\frac{\beta_{t-1}}{t-1}\right)\frac{\beta_t}{t}.$$
We verify easily that 
${\mathbb E}\, \mu_t =\gamma_t\mu_0+(1-\gamma_t)\bar\mu$, for $t>0$;
therefore, the first part in establishing $\eps$-self-sustainability 
consists of proving that  
\begin{equation}\label{gamma_t<>}
1\geq\gamma_t\geq 1- \frac{\eps}{ |\mu_0-\bar\mu|}\, ,
\end{equation}
which will show that $|{\mathbb E}\, \mu_t - \mu_0|\leq \eps$.
Note that 
$$\gamma_t\leq  \frac{1}{t}
\prod_{s=1}^{t-1}\Bigl( 1 + \frac{1}{s}\Bigr) =1.
$$
Now define
$$
\alpha_s= \frac{\eps}{  B_c |\mu_0-\bar\mu| s (1+\log s)^{1+c}}.
$$
for $s>0$. 
We pick a constant $B_c$ large enough so that 
$\alpha_s$ is small enough to carry out first-order Taylor approximations 
around $1+\alpha_s$. We find that
\begin{equation*}
\begin{split}
1+\frac{\beta_s}{s}
&= 1+\frac{1}{s}\Bigl(1-\frac{1}{1+m_t\tau/\bar\tau}\Bigr)
\geq \Bigl(1+\frac{1}{s}\Bigr)\Bigl( 1- \frac{1}{(s+1)m_t\tau/\bar\tau}\Bigr) \\
&\geq \Bigl(1+\frac{1}{s}\Bigr)
\Bigl(1- \frac{s \alpha_s}{s+1}\Bigr)
\geq \Bigl(1+\frac{1}{s}\Bigr)(1-\alpha_s)
\geq \Bigl(1+\frac{1}{s}\Bigr)e^{-2 \alpha_s}.
\end{split}
\end{equation*}
Thus,
\begin{equation*}
\gamma_t
\geq \frac{\beta_t}{t}
\prod_{s=1}^{t-1} \Bigl(1 + \frac{1}{s} \Bigr)
     e^{-2\sum_{s=1}^{t-1} \alpha_s}
= \beta_t  e^{-2\sum_{s=1}^{t-1} \alpha_s} 
\geq 1 - \frac{\eps}{|\mu_0-\bar\mu|},
\end{equation*}
which establishes~(\ref{gamma_t<>}).
Our derivation relies on the fact that
$$
 \beta_t \geq  1 - \frac{\eps}{  B_c |\mu_0-\bar\mu|  (1+\log t)^{1+c}}
 \geq 1 - \frac{\eps}{  2 |\mu_0-\bar\mu| }
$$
and
\begin{equation*}
\sum_{s=1}^{t-1} \frac{1}{s(1+\log s)^{1+c}}
\leq 1+ \frac{1}{(\log e)^{1+c}} \int_2^{t-1} \frac{1}{x(\ln x)^{1+c}}\, dx 
= O\Bigl( \frac{1}{c} \Bigr)\, ;
\end{equation*}
hence,
$$  e^{-2\sum_{s=1}^{t-1} \alpha_s} 
\geq  
 e^{- O(\eps / ( cB_c |\mu_0-\bar\mu| ))}
 \geq 1 - \frac{\eps}{  2 |\mu_0-\bar\mu| }.
$$

Having shown that $|{\mathbb E}\, \mu_t - \mu_0|\leq \eps$ for all $t$,
it now suffices to prove that $\sigma_t^2 + \Varr \mu_t$ remains bounded. 
We note that $\tau_t>m_t\tau\rightarrow\infty$, hence $\sigma_t^2=1/\tau_t\rightarrow 0$, so
the remainder of the proof needs to establish that the variance of $\mu_t$ stays bounded.
Writing $D_{t,s}:=   d_{t,s,1}+\dots+d_{t,s,m_t}$, we have
$\Varr  D_{t,s} =m_t\Varr  d_{t,s,1} =m_t(\sigma_s^2+\sigma^2+\Varr  \mu_s )$;
hence
$${\mathbb E}\, D_{t,s}^2 =\Varr  D_{t,s} +({\mathbb E}\, D_{t,s})^2=
m_t(\sigma_s^2+\sigma^2+\Varr  \mu_s )+m_t^2({\mathbb E}\, \mu_s)^2.
$$
In~(\ref{random_update}), 
the variables $\chi_{t,s}$ and $D_{t,s}$ are independent,
for $0\leq s\le t-1$; furthermore,
${\mathbb E}\, \chi_{t,s}={\mathbb E}\, \chi_{t,s}^2=1/t$, and ${\mathbb E}\, \chi_{t,s_1}\chi_{t,s_2}=0$ if $s_1\neq s_2$;
therefore,   
\begin{equation}\label{varIneq}
\begin{split}
\Varr [ \chi_{t,s} D_{t,s}] 
&=
{\mathbb E}\, \chi_{t,s}^2\, {\mathbb E}\, D_{t,s}^2 
- ({\mathbb E}\, \chi_{t,s})^2    ({\mathbb E}\,  D_{t,s})^2 
=   \frac{{\mathbb E}\, D_{t,s}^2}{t}    - \frac{ ({\mathbb E}\, D_{t,s})^2  }{t^2} \\
&=
\Bigl(\frac{ m_t}{t}\Bigr) \bigl(\sigma_s^2+\sigma^2+\Varr  \mu_s 
+m_t ({\mathbb E}\, \mu_s)^2 \bigr)
-\Bigl(\frac{ m_t}{t}\Bigr)^2 ({\mathbb E}\, \mu_s)^2  
\end{split}
\end{equation}
and, for $s_1\neq s_2$,
\begin{equation}\label{covIneq}
\begin{split}
\Cov [\chi_{t,s_1}D_{t,s_1},\chi_{t,s_2}D_{t,s_2}]
&=\,\, {\mathbb E}\, [\chi_{t,s_1}\chi_{t,s_2}D_{t,s_1},D_{t,s_2}]
   -{\mathbb E}\, [\chi_{t,s_1}D_{t,s_1}]{\mathbb E}\, [\chi_{t,s_2}D_{t,s_2}]\\
&=\,\,
{\mathbb E}\, [\chi_{t,s_1}\chi_{t,s_2}]{\mathbb E}\, [D_{t,s_1}D_{t,s_2}]-{\mathbb E}\, \chi_{t,s_1}{\mathbb E}\, D_{t,s_1}{\mathbb E}\, \chi_{t,s_2}{\mathbb E}\, D_{t,s_2}\\
&= \,\, -\frac{1}{t^2}\, {\mathbb E}\, D_{t,s_1}{\mathbb E}\, D_{t,s_2}
= -\Bigl( \frac{m_t}{t}\Bigr)^2 {\mathbb E}\, \mu_{s_1} {\mathbb E}\, \mu_{s_2}.
\end{split}
\end{equation}
Then, by taking the variance on both sides of \eqref{random_update}, we have
\begin{equation*}
\begin{split}
\Varr \mu_t 
&=\Bigl( \frac{\beta_t}{m_t}\Bigr)^2 \, \Varr \sum_{s=0}^{t-1}\chi_{t,s}D_{t,s} \\
&=\Bigl( \frac{\beta_t}{m_t}\Bigr)^2 \,
\biggl( \sum_{s=0}^{t-1}\Varr  [\chi_{t,s}D_{t,s}]+\sum_{0\le s_1\neq s_2\le t-1}\Cov [\chi_{t,s_1}D_{t,s_1},\chi_{t,s_2}D_{t,s_2}] \biggr) \\
&=\Bigl( \frac{\beta_t}{m_t}\Bigr)^2 \,
\biggl(\sum_{s=0}^{t-1}\Bigl(\frac{ m_t}{t}\Bigr) \bigl(\sigma_s^2+\sigma^2+\Varr  \mu_s 
+m_t ({\mathbb E}\, \mu_s)^2 \bigr)
-\Bigl(\frac{ m_t}{t}\Bigr)^2 \Bigl(\sum_{s=0}^{t-1} {\mathbb E}\, \mu_s\Bigr)^2
\biggr) \\
&\leq \frac{1}{t m_t}
      \sum_{s=0}^{t-1}
          \bigl(\sigma_s^2+\sigma^2+\Varr  \mu_s 
                     +m_t ({\mathbb E}\, \mu_s)^2 \bigr) .
\end{split}
\end{equation*}
Notice that $\sigma_s^2\rightarrow 0$ and $({\mathbb E}\, \mu_s)^2$ 
is bounded since $|{\mathbb E}\, \mu_t-\mu_0|\le  \eps $. 
We conclude that $\sigma_t^2 + \Varr \mu_t$ remains bounded for all $t$.
\hfill $\Box$
\proofend

\section{Iterated Bayesian Linear Regression}\label{sec:linear}

The iterated version of Bayesian linear regression has been 
the subject of extensive study in the field of psychology~\cite{kalish2007iterated,beppu2009iterated,tamariz2015culture,
bartlettremembering,griffiths2008theoretical}.
The work has involved experimentation with human subjects but little in the way 
of theoretical analysis. This section is a first step toward filling this void. 
The task at hand is to estimate a hypothesis $h\in {\mathcal H}:={\mathbb R}^d$
given a noisy measurements on the hyperplane $y=h^T x$, where $x\in {\mathbb R}^d$. 
In the Bayesian setting, we assume a Gaussian prior on the hypothesis space:
${\mathbb P}[h]\sim N(\bar{\mu},\bar{\sigma}^2 I_d)$.
The data is given by $(x,y)$, where $x\sim N(0,I_d)$ and
$y= h^T x+\phi$, for $\phi\sim N(0,\sigma^2)$ (with $x,\phi$ independent).
Since we typically make several measurements, we write this (likelihood) relation in matrix form:
$ y= Xh + \phi$, where $y\in {\mathbb R}^m$ (with $m$ the number of measurements);
$\phi\sim  N(0, \sigma^2 I_m)$; and $X$ is an $m$-by-$d$ matrix each of whose rows denotes 
a random vector $x\sim N(0,I_d)$.  This means that the matrix $X$ is random (a fact of key
importance in our discussion below). We have:

\begin{equation*}
\begin{cases}
\,\, {\mathbb P}[\phi]\sim \,  \exp \bigl \{ -\frac{1}{2\sigma^2}\| \phi\|_2^2 \bigr\}  
\hspace{4.2cm} \text{{\small\tt (noise) }}
       \vspace{.2cm} \\
\,\, {\mathbb P}[h] \sim \,     \exp \bigl \{ -\frac{1}{2\bar \sigma^2} \|h - \bar\mu\|_2^2 \bigr\}  
\hspace{3.5cm} \text{{\small\tt (prior) }}
           \vspace{.2cm}  \\
\,\, {\mathbb P}[y |X,h] \sim \,    \exp  \bigl \{ -\frac{1}{2 \sigma^2} \|y- Xh\|_2^2 \bigr\} 
\hspace{1.5cm} \text{{\small\tt (likelihood) }} 
\end{cases}
\end{equation*}

\bigskip

In iterated Bayesian linear regression, the $t$-th learner receives 
her data from learner $t-1$. Here, learner $0$ is treated just like any other
agent, except that his prior ${\mathbb P}[h]\sim N(\mu_0,\bar{\sigma}^2 I_d)$
is the distribution to be learned iteratively.
Since sampling from the prior is independent of $X$, 
Bayesian updating gives the posterior $N(\mu_t, \Sigma_t)$, where
\begin{equation*}
{\mathbb P}[h|X,y]=  {\mathbb P}[h] \,{\mathbb P}[y|X,h]/  {\mathbb P}[y|X]
\sim 
  \exp \Bigl \{ -\frac{1}{2\bar \sigma^2} \|h - \bar\mu\|_2^2 -
              \frac{1}{2 \sigma^2} \|y- Xh\|_2^2 \Bigr\} .
\end{equation*}
Completing the square in the usual fashion shows that 
the posterior of learner $t$ is given by:
\begin{equation}\label{linearUpdate1}
\begin{cases}
\,\, \Sigma_{t} =
     \left( \bar{\sigma}^{-2} I_d + \sigma^{-2} X_t^TX_t \right)^{-1}
                   ;  \vspace{.2cm} \\
\,\, \mu_{t} = \Sigma_{t}\left ( \bar{\sigma}^{-2}\bar\mu+ \sigma^{-2}X_t^Ty_t \right),
\end{cases}
\end{equation}
where $(X_t,y_t)$ is the data gathered by learner $t$ from her predecessor:
specifically, $y_t=X_t h+\phi_t$,  where
$h$ is collected from the $(t-1)$-th learner by sampling his posterior distribution 
$N(\mu_{t-1},\Sigma_{t-1})$.

\bigskip
\begin{theorem}\label{iterLR}
Given any small enough $\delta, \eps>0$, 
the following sample size sequence for iterated Bayesian linear regression ensures that  
$\|{\mathbb E}\mu_t-\mu_0\|_2\le  \delta $
with probability greater than $1-\eps$:
$$m_t=    D_c \frac{ \|\mu_0-\bar{\mu}\|_2}{\delta} 
                   \Bigl( \frac{\sigma}{\bar{\sigma}} \Bigr)^{\!2}  
                        t^{1+c} 
                        +D_c \,  d\log \frac{t+1}{\eps},
$$
for an arbitrarily small $c>0$ and a constant $D_c$ that depends only on c.
\end{theorem}
\proof
We proceed in two steps: first, we show that to keep ${\mathbb E}\mu_t$ arbitrarily
close to $\mu_0$ for all $t$ hinges on spectral properties of certain random
matrices; second, we call on known facts about the singular values of
random Gaussian matrices to translate the spectral condition into a high-probability event.
The proof unfolds as a series of simple relations, which we state first and then demonstrate.
The first one follows directly from~\eqref{linearUpdate1}:

\begin{equation}\label{linear_update}
{\mathbb E}\mu_{t}=\left(I_d+M_t\right)^{-1} \left(\bar{\mu}+M_t \, {\mathbb E}\mu_{t-1}\right),
\hspace{.2cm}
\text{where}
\hspace{.5cm}
M_t:=\left(\frac{\bar{\sigma}}{\sigma}\right)^2X_t^TX_t.
\end{equation}
Note that \eqref{linear_update} is 
a randomized recursive relation since the data points $X_1,X_2,\dots$ are themselves random.
We note that all the matrices whose inverses are taken are positive definite, hence nonsingular.
To move on to our second relation, we define the matrix
$$Q_t:=(I_d+M_{t})^{-1}M_{t} (I_d+M_{t-1})^{-1}M_{t-1}\cdots(I_d+M_1)^{-1}M_1, $$
for $t>0$, with $Q_0= I_d$, 
and prove by induction that
\begin{equation}\label{EQt}
{\mathbb E}\, \mu_t =Q_t\mu_0+(I_d-Q_t)\bar{\mu}.
\end{equation}
The base case is obvious so we assume that $t>0$: 
by~(\ref{linear_update}),
\begin{equation*}
\begin{split}
{\mathbb E}\, \mu_t 
&=  (I_d + M_t)^{-1}(\bar\mu + M_t\, {\mathbb E}\, \mu_{t-1})  \\
&=  (I_d + M_t)^{-1}(\bar\mu + M_tQ_{t-1}\mu_0 + M_t(I_d- Q_{t-1})\bar\mu )\\
&=  (I_d + M_t)^{-1} M_t Q_{t-1}\mu_0 
+ (I_d + M_t)^{-1}(I_d+ M_t (I_d-Q_{t-1}))\bar\mu \\
&= Q_t \, \mu_0 + ( I_d- (I_d+M_t)^{-1}M_t Q_{t-1})\bar\mu,
\end{split}
\end{equation*}
which proves~(\ref{EQt}).
Our next goal is to bound the information decay
$\|{\mathbb E}\, \mu_t -\mu_0\|_2$. To do that, we investigate 
the spectral norm of the matrix $I_d- Q_t$, which  leads to our
third relation. We prove by induction that, for $t>0$,
\begin{equation}\label{spectralBnd}
\|I_d- Q_t\|_2\leq \sum_{s=1}^{t} \|A_s\|_2,
\end{equation}
where $A_s:= (I_d+ M_s)^{-1}$.
For $t=1$, $Q_1= (I_d+M_1)^{-1}M_1= I_d- (I_d+M_1)^{-1}$ and the
claim follows. If $t>1$, then
\begin{equation*}
\begin{split}
\|I_d- Q_t\|_2
&= \| (I_d- Q_{t-1})+  (Q_{t-1}-Q_t) \|_2 \\
&\leq \| I_d- Q_{t-1}\|_2 + \| Q_{t}-Q_{t-1} \|_2 
\leq 
\sum_{s=1}^{t-1} \|  A_s \|_2 +
\|  \Psi \|_2,
\end{split}
\end{equation*}
where 
$\Psi:= (A_tM_t- I_d)Q_{t-1}$.
Since $A_{t}(I_d+ M_{t}) = I_d$, we have
$\Psi= - A_{t}Q_{t-1}$.
Each matrix $M_s$ is positive semidefinite, so the eigenvalues of
$(I_d+ M_s)^{-1}M_s$ are of the form $\lambda/(1+\lambda)$, where $\lambda\geq 0$.
This shows that all the eigenvalues of $Q_s$ are between 0 and 1; therefore
$\|Q_s\|_2\leq 1$. The eigenvalues of $I_d- A_tM_t$ are the same as those of $A_t$; hence,
by submultiplicativity,
$\|\Psi\|_2\leq \|A_t\|_2 \|Q_{t-1}\|_2\leq \|A_t\|_2$,
which establishes~(\ref{spectralBnd}). 

We are now ready to express the information decay in spectral terms.
Pick an arbitrarily small constant $c>0$ and assume that
\begin{equation}\label{CondIM}
\|A_s\|_2\leq \frac{\delta}{\|\bar\mu - \mu_0\|_2}
\Bigl(\frac{c}{1+c}\Bigr)
\Bigl(\frac{1}{s}\Bigr)^{1+c}.
\end{equation}
By~(\ref{EQt}), ${\mathbb E}\, \mu_t - \mu_0 = (I_d- Q_t)(\bar\mu - \mu_0)$; therefore,
by~(\ref{spectralBnd}),
\begin{equation}\label{Emu-mu0}
\begin{split}
\|{\mathbb E}\, \mu_t -\mu_0\|_2
&\leq  \|\bar{\mu} - \mu_0\|_2  \sum_{s=1}^{t} \|A_s\|_2 
\leq \frac{\delta c}{1+c} \sum_{s=1}^{t} s^{-1-c} \\
&\leq  \frac{\delta c}{1+c}\Bigl( 1 + \int_1^\infty x^{-1-c}\, dx \Bigr)
%%% =  \frac{\delta c}{1+c}\Bigl(1+ \frac{1}{c}\Bigr)
= \delta,
\end{split}
\end{equation}

\noindent
The relation says that, on average, the means of 
any of the agents' posteriors
can be brought as close to the original mean to be learned as we want.
We can turn this into a high-probability event by using some basic
random matrix theory.  Recall that 
${\mathbb E}\, \mu_t$ is itself a random variable whose stochasticity
comes from the matrices $X_s$, which are all drawn from Gaussians.
Because $M_s$ is positive semidefinite, 
\begin{equation}\label{AsUB}
\|A_s\|_2\le\|M_s^{-1}\|_2\leq
    \frac{ (\sigma/\bar{\sigma})^2}{ \lambda_{\text{min}} (X_t^TX_t)}
    \leq \Bigl( \frac{\sigma/\bar{\sigma}}{\sigma_1(X_t)} \Bigr)^2,
\end{equation}
which gives us a relation between
the spectral norm of $(I_s+M_s)^{-1}$ and 
the smallest singular value $\sigma_1(X_t)$ of an $m_t$-by-$d$ matrix $X_t$
whose elements are drawn {\em iid} from $N(0,1)$.
The asymptotic behavior of $\sigma_1(X_t)$ for large values of $m_t$ 
has been extensively studied within the field of random matrix theory~\cite{davidson2001local,edelman1988eigenvalues,rudelson2009smallest}.
Following Theorem II.13 in (Davidson \& Szarek~\cite{davidson2001local}),
for any $\gamma_t>0$, 
$${\mathbb P}[\sigma_1(X_t)<\sqrt{m_t}-\sqrt{d}-\gamma_t]\le e^{-\gamma_t^2/2}.$$
We use $C$ below as a generic constant large enough to satisfy the inequalities where it appears.
Setting $\gamma_t = C\sqrt{ \log ((t+1)/\eps) }$ ensures that 
$\sum_{t>0} e^{-\gamma_t^2/2} <\eps$, hence that
$\sigma_1(X_t)<\sqrt{m_t}-\sqrt{d}-\gamma_t$ holds for all $t$ with probability 
less than $\eps$. With our setting of $m_t$, this means that, for all $t>0$,
\begin{equation}\label{HighProbSigma1}
{\mathbb P}\Bigl[\, \sigma_1(X_t)\geq  \frac{\sqrt{m_t}}{2} \, \Bigr] > 1-\eps.
\end{equation}
Assuming the event in~(\ref{HighProbSigma1}),
it follows from~(\ref{AsUB}) and our setting of $m_t$ that
$$
\|A_t\|_2\leq  \frac{4}{m_t}  \Bigl(\frac{\sigma}{\bar{\sigma}}\Bigr)^2
\leq \frac{\delta}{\|\bar\mu - \mu_0\|_2}
\Bigl(\frac{4}{D_c}\Bigr)
\Bigl(\frac{1}{t}\Bigr)^{1+c};
$$
hence~(\ref{CondIM}) for $D_c$ large enough.
By~(\ref{Emu-mu0}, \ref{HighProbSigma1}), this proves that, with probability greater than $1-\eps$,
$\|{\mathbb E}\, \mu_t -\mu_0\|_2 \leq \delta$ for all $t>0$,
which completes the proof.
\hfill $\Box$
\proofend

\newpage

\bibliography{refer}
\bibliographystyle{unsrt}

\newpage

\begin{center}
\huge\bf Appendix
\end{center}

\vspace{1cm}

The two forms of the function $d_{RS}$ in~(\ref{d_RS}) make it clear that 
$0\le d_{RS}(\bm{a},\bm{b}) \le 1$ and $d_{RS}(\bm{a},\bm{b})=0$ 
if and only if $\bm{a}$ and $\bm{b}$ are identical.
We easily check that $d_{RS}$
makes the simplex ${\mathcal S}$ of distributions over ${\mathcal D}$ into a metric space.
Indeed, $d_{RS}(\cdot,\cdot)$ is obviously symmetric, and $d_{RS}(\bm{a},\bm{b})=0$ 
implies that $\bm{a}=\bm{b}$.
To check the triangular inequality, notice that
\begin{equation}\label{angle}
d_{RS}(\bm{a},\bm{b})=\sqrt{1-\Bigl( \sum_{i=1}^s \sqrt{ a_ib_i}\Bigr)^2}=\sin \langle \sqrt{\bm{a}},\sqrt{\bm{b}}\,\rangle,
\end{equation}
where $\langle \sqrt{\bm{a}},\sqrt{\bm{b}}\,\rangle$ is the angle between the unit vectors 
$\sqrt{\bm{a}}$ and $\sqrt{\bm{b}}$, using the notation
$\sqrt{\bm{v}}=(\sqrt{v_1},\dots,\sqrt{v_s})$.
To prove that $d_{RS}(\bm{a},\bm{b})+d_{RS}(\bm{b},\bm{c})\ge d_{RS}(\bm{a},\bm{c})$ for any $\bm{a},\bm{b},\bm{c}\in{\mathcal S}$,
we denote by $\alpha, \beta, \gamma$ the corresponding angles in that order,
ie,  $\alpha=\langle \sqrt{\bm{a}},\sqrt{\bm{b}}\,\rangle$, etc.
The coordinates in $\bm{a},\bm{b},\bm{c}$ are nonnegative; therefore $0\le \alpha,\beta,\gamma\le \pi/2$.
These form the three angles at the origin of a tetrahedron with a vertex at the origin;
therefore, by the triangular inequality in spherical geometry, $\alpha+\beta\ge\gamma$.
If $\alpha+\beta\le \frac{\pi}{2}$, then $\sin\alpha+\sin\beta\ge \sin\alpha\cos\beta+\cos\alpha\sin\beta=\sin(\alpha+\beta)\ge\sin\gamma$.
On the other hand, if $\alpha+\beta>\pi/2$, then 
$\sin\alpha+\sin\beta=2\sin\frac{\alpha+\beta}{2}\cos\frac{\alpha-\beta}{2}\ge 2\sin\frac{\pi}{4}\cos\frac{\pi}{4}=1\ge \sin\gamma$, which establishes
the triangular inequality.

\paragraph{Relation to the Euclidean distance.}\
Shrinking the simplex ${\mathcal S}$ by a tiny amount, we define
${\mathcal S}_ \eps := \{\bm{a}\in{\mathcal S}:  \eps \le a_i\le 1- \eps \}$
and note that
$$d_{E}(\bm{a},\bm{b}):=\|\bm{a}-\bm{b}\|_2
= \sqrt{ \sum_{i=1}^s(\sqrt{a_i}-\sqrt{b_i})^2(\sqrt{a_i}+\sqrt{b_i})^2 }.
$$
It follows that, for $\bm{a},\bm{b}\in {\mathcal S}_ \eps $,

\begin{equation}\label{DistRelat}
\frac{1}{2} d_E(\bm{a},\bm{b}) \le d_E(\sqrt{\bm{a}},\sqrt{\bm{b}}\,)\le \frac{1}{2\sqrt{ \eps }}\,d_E(\bm{a},\bm{b}).
\end{equation}
On the other hand,
$\|\sqrt{\bm{a}}\|_2=\|\sqrt{\bm{b}}\|_2=1$,
so the vectors $\sqrt{\bm{a}}$ and $\sqrt{\bm{b}}$ form an isosceles triangle; hence
$$d_E(\sqrt{\bm{a}},\sqrt{\bm{b}}\,)=2\sin \frac{1}{2}\langle\sqrt{\bm{a}},\sqrt{\bm{b}}\rangle=\frac{\sin\langle\sqrt{\bm{a}},\sqrt{\bm{b}}\rangle}{\cos \frac{1}{2}\langle\sqrt{\bm{a}},\sqrt{\bm{b}}\rangle}
=\frac{d_{RS}(\bm{a},\bm{b})}{\cos\frac{1}{2}\langle\sqrt{\bm{a}},\sqrt{\bm{b}}\rangle}\, .$$
Since $0\le \langle\sqrt{\bm{a}},\sqrt{\bm{b}}\,\rangle\le\frac{\pi}{2}$,
$$d_{RS}(\bm{a},\bm{b})\le d_E(\sqrt{\bm{a}},\sqrt{\bm{b}}\,)\le \sqrt{2}\, d_{RS}(\bm{a},\bm{b}).$$
Together with~(\ref{DistRelat}) this shows that, for any
$\bm{a},\bm{b}\in {\mathcal S}_ \eps $,
\begin{equation}\label{dE-dRS}
\frac{1}{2\sqrt{2}}\, d_E(\bm{a},\bm{b})\le d_{RS}(\bm{a},\bm{b})  
\le \frac{1}{2\sqrt{ \eps }}\, d_E(\bm{a},\bm{b}),
\end{equation}
which shows that the Euclidean distance and the metric $d_{RS}$ are equivalent in ${\mathcal S}_ \eps $.

\paragraph{Relation to other distances.}\
The metric $d_{RS}$ is related to the Hellinger and Bhattacharyya distances.
Writing $C(\bm{a},\bm{b})=\sum_{i=1}^s\sqrt{a_ib_i}$~\cite{comaniciu2000real}, 
then $d_{RS}(\bm{a},\bm{b})=\sqrt{1-C(\bm{a},\bm{b})^2}$.
The Hellinger distance is defined as $d_H(\bm{a},\bm{b})=\sqrt{1-C(\bm{a},\bm{b})}$~\cite{hazewinkel2013encyclopaedia}, while
the Bhattacharyya distance is defined as $d_B(\bm{a},\bm{b})=-\ln C(\bm{a},\bm{b})$~\cite{bhattachayya1943measure}.
The total variation distance $d_{TV}$ is
half the $\ell_1$-norm; therefore
$d_{TV}(\bm{a},\bm{b})\leq \frac{1}{2}\sqrt{s}\, d_E(\bm{a},\bm{b})$.
Combining these observations with~(\ref{dE-dRS}) establishes~(\ref{DistCompare}):

\begin{equation*}
\begin{cases}
\, d_H= \sqrt{1-\sqrt{1-d_{RS}^2}}\, ; \vspace{.2cm} \\
\, d_B= -\frac{1}{2}\ln (1- d_{RS}^2)\, ; \vspace{.2cm} \\
\, d_{TV}\leq \sqrt{2s} \, d_{RS}.
\end{cases}
\end{equation*}

\end{document}